\documentclass[a4paper,12pt]{article}

\makeatletter
\@addtoreset{footnote}{page}
\makeatother

\renewenvironment{enumerate}
{\begin{list}{\parbox{2em}{(\arabic{enumi})}}{
\usecounter{enumi}
\setlength{\topsep}{0ex}
\setlength{\itemindent}{1em}
\setlength{\leftmargin}{2em} %left indent
\setlength{\rightmargin}{1em}%right indent
\setlength{\labelsep}{1em}%sep between label and text
\setlength{\labelwidth}{3em}%
\setlength{\itemsep}{0.1em}%line break
\setlength{\parsep}{0em}%paragraph break
\setlength{\listparindent}{0em}%paragraph indent
}}{\end{list}}

\usepackage{amsthm}
\usepackage{amsmath,amssymb,latexsym,amsfonts,mathrsfs}

\renewcommand{\d}{{\rm d}} %differential
\newcommand{\Var}{\mathop{\rm Var}}
\newcommand{\Cov}{\mathop{\rm Cov}}

\newcommand{\dist}{\stackrel{{\rm d}}{=}}

\newcommand{\tend}[2]{\mathrel{\mathop{\longrightarrow}\limits^{#1}_{#2}}}

\newcommand{\supre}[2]{\mathrel{\mathop{\sup}\limits^{#1}_{#2}}}

\newcommand{\absol}[1]{\left| #1 \right|} %absolute value
\newcommand{\norm}[1]{\left\| #1 \right\|} %norm
\newcommand{\rbra}[1]{\!\left( #1 \right)} %round brackets or parentheses
\newcommand{\cbra}[1]{\!\left\{ #1 \right\}} %curly brackets or braces
\newcommand{\sbra}[1]{\!\left[ #1 \right]} %square brackets or brackets
\newcommand{\bE}{\ensuremath{\mathbb{E}}}

\newcommand{\bP}{\ensuremath{\mathbb{P}}}

\newcommand{\cF}{\ensuremath{\mathcal{F}}}

\newcommand{\cP}{\ensuremath{\mathcal{P}}}

\theoremstyle{plain}
\newtheorem{Thm}{Theorem}[section]

\newtheorem{Lem}[Thm]{Lemma}

\newtheorem{Cor}[Thm]{Corollary}

\theoremstyle{definition}

\newtheorem{Def}[Thm]{Definition}

\newcommand{\Proof}[2][Proof]{\begin{proof}[{#1}] #2 \end{proof}}

\numberwithin{equation}{section}

\makeatletter
\renewcommand\section{\@startsection {section}{1}{\z@}%
                                   {-3.5ex \@plus -1ex \@minus -.2ex}%
                                   {2.3ex \@plus.2ex}%
                                   {\normalfont\large\bf}}
\makeatother

\makeatletter
\renewcommand\subsection{\@startsection {subsection}{1}{\z@}%
                                   {-3.5ex \@plus -1ex \@minus -.2ex}%
                                   {2.3ex \@plus.2ex}%
                                   {\normalfont\normalsize\bf}}
\makeatother
\usepackage{graphicx}  % For including images
\usepackage{caption}   % For customizing captions
\usepackage{float}     % For precise placement of figures

\begin{document}

\begin{center}
{\Large \bf 
Law of large numbers and central limit theorem for renewal Hawkes processes
}
\end{center}
\footnotetext{
This research was supported by RIMS and by ISM.}
\begin{center}
Luis Iv\'an Hern\'andez Ru\'{\i}z \footnote{
Graduate School of Science, Kyoto University. Contact: luisivanhr@ciencias.unam.mx}\footnote{
The research of this author was supported
by JSPS Open Partnership Joint Research Projects grant no. JPJSBP120209921.
}
\end{center}
\begin{abstract}
A uniform law of large numbers and a central limit theorem are established via martingale and cluster process techniques for a univariate Hawkes process with immigration given by a renewal process. The results are obtained for renewal processes with absolutely continuous interarrival distribution.
\end{abstract}

\section{Introduction}
Hawkes \cite{HawkesOriginal} introduced \emph{classical Hawkes processes} as point processes with a self-exciting nature, in the sense that previous events facilitate the ocurrence of future events. Hawkes--Oakes \cite{HawkOak} showed that the classical Hawkes process could be understood as an \emph{independent cluster process} in which the centre process is given by a homogeneous Poisson process of immigrants and the satellite processes are given by branching processes formed by the offspring of those immigrants. As a generalization to the classical case, Wheatley--Filimonov--Sornette \cite{Wheatley2016} introduced the \emph{Renewal Hawkes process} (abbreviated RHP) in which immigration is given by a renewal process. This generalization allows for more flexibility when fitting Hawkes processes to data sets as Stindl--Chen \cite{StindlWeibull} did for modelling financial returns using Hawkes processes where the renewals were Weibull distributed. Chen--Stindl \cite{ChenDirect} and Stindl--Chen \cite{StindlMulti}  studied  the evaluation of the likelihood for the single variate and multivariate RHP respectively, and explained the challenges of computing the likelihood with respect to the natural filtration, and Chen--Stindl \cite{ChenImprove} refined the method of evaluation to improve the speed of the calculation. Hern\'andez--Yano \cite{IvanYanoConstruction2023} established a representation of the RHP as an independent cluster process and computed the  \emph{probability generating functional} for the RHP and its stationary limit process.

The purpose of this paper is to establish limit theorems for the RHP, namely, a law of large numbers and a central limit theorem. For the law of large numbers we use a martingale approach similar to what Bacry--Delattre--Hoffmann--Muzy \cite{BacryLimit} did for the classical Hawkes process. In the case of the central limit theorem, we make use of the cluster representation of the RHP described by Hern\'andez--Yano \cite{IvanYanoConstruction2023}. We tackle the univariate case and make some comments at the end of the paper about the possibility of generalizing it to the multivariate case.

Let $\cbra{T_n}_{n\ge0}$ be a sequence of random variables on $[0,\infty)$ defined on a common probability space $(\Omega,\cF,\bP)$ and such that $T_0=0$ and, for all $i\ge0$, we have $T_i<T_{i+1}$ on the event $\cbra{T_i<\infty}$ and $T_{i+1}=\infty$ on $\cbra{T_i=\infty}$. We identify the \emph{point process} $\cbra{T_n}_{n\ge0}$ with the associated counting process $N(t)=\sum_{i}1_{\cbra{T_i\le t}}$ for $t\ge0$. Let $(\cF_t)_{t\ge0}$ be a filtration to which $N$ is adapted. We can specify $N$ through its \emph{intensity}.

\begin{Def}
\label{DefIntensity}
Let $N$ be a point process and $(\cF_t)_{t\ge0}$ a filtration to which $N$ is adapted. Let $\lambda$ be a nonnegative, a.s. locally integrable process that is $(\cF_t)$-progressive. We say that $N$ admits the $(\cF_t)$-\emph{intensity} $\lambda$ if the process given as
\begin{align*}
M(t)=N(t)-\int_0^t\lambda(s)\d s,\quad t\ge0,
\end{align*}
is an $(\cF_t)$-martingale, in which case $M$ is called the \emph{characteristic martingale} of $N$. 
\end{Def}
As a consequence of Definition \ref{DefIntensity}, for any nonnegative process $C(t)$ that is \emph{predictable}, i.e. for all $t\ge0$ it is measurable with respect to the $\sigma$-field 
\begin{align}
\cP\rbra{\cF_t}=\sigma\rbra{(s,t]\times A; 0\le s\le t, A\in\cF_s},
\end{align}
it holds that
\begin{align}
\label{PropertyIntensity}
\bE\sbra{\int_0^\infty C(s)N(\d s)}=\bE\sbra{\int_0^\infty C(s)\lambda(s)\d s}.
\end{align}
\subsection{Limit theorems for the classical Hawkes process}
The classical Hawkes process is defined through its intensity.
\begin{Def}A point process $N$ is called a \emph{classical (univariate) Hawkes process} if $N$ admits an $\rbra{\cF_t}$-intensity given as
\begin{align}
\label{classicalHawkes}
\lambda(t)=\mu+\int_0^{t}h(t-u)N(\d u),\quad t\ge0,
\end{align}
where $\mu$ is a positive constant and $h$ is a nonnegative measurable function on $[0,\infty)$ satisfying $\int_0^\infty h(t)\d t<1$.
\end{Def}
Bacry--Delattre--Hoffmann--Muzy \cite{BacryLimit} found limit theorems for classical multivariate Hawkes processes. In the case of one dimension, if $\alpha:=\int_0^\infty h(t)\d t<1$, a law of large numbers (LLN) is given as
\begin{Thm}[\textbf{Bacry--Delattre--Hoffmann--Muzy \cite[Theorem 1]{BacryLimit}}]
\label{BacryLimT1}
We have \hfill \break$N(t)\in L^2(\bP)$ for all $t\ge0$ and the convergence
\begin{align}
\supre{}{v\in[0,1]}\absol{T^{-1}N(Tv)-v\frac{\mu}{1-\alpha}}\tend{}{T\rightarrow\infty}0
\end{align}
holds a.s. and in $L^2(\bP)$.
\end{Thm}
A central limit theorem (CLT) was also proved, namely,
\begin{Thm}[\textbf{Bacry--Delattre--Hoffmann--Muzy \cite[Theorem 2]{BacryLimit}}]
\label{BacryLimT2}
The convergence
\begin{align}
\rbra{\frac{1}{\sqrt{T}}\rbra{N_{Tv}-\bE\sbra{N_{Tv}}}}_{v\in[0,1]}\tend{d}{T\rightarrow\infty}\rbra{\sqrt{\frac{\mu}{(1-\alpha)^3}}W_v}_{v\in[0,1]}
\end{align}
holds in the Skorokhod topology, where $(W_v)_{v\in[0,1]}$ is a standard Brownian motion.
\end{Thm}
And finally, assuming that $\int_0^\infty t^{1/2}h(t)\d t<\infty$, the result below for asymptotic normality also follows.
\begin{Cor}[\textbf{Bacry--Delattre--Hoffmann--Muzy \cite[Corollary 1]{BacryLimit}}]
\label{CorBacry}
The convergence
\begin{align}
\rbra{\frac{1}{\sqrt{T}}N(Tv)-v\frac{m}{1-\alpha}\sqrt{T}}_{v\in[0,1]}\tend{d}{T\rightarrow\infty}\rbra{\sqrt{\frac{\mu}{(1-\alpha)^3}} W(v)}_{v\in[0,1]},
\end{align}
holds in the Skorokhod topology, where $(W_v)_{v\in[0,1]}$ is a standard Brownian motion.
\end{Cor}
\subsection{Limit theorems for the renewal Hawkes process}
We extend these results to the RHP. This process has a cluster structure that has been discussed in Hern\'andez--Yano \cite{IvanYanoConstruction2023} and which we briefly review here. Consider the following assumptions
\begin{enumerate}
\item[\textbf{(A0)}]$h$ is a nonnegative measurable function on $[0,\infty)$ satisfying $\alpha:=\int_0^\infty h(t)\d t<1$.
\item[\textbf{(B0)}]$F$ is a probability distribution on $[0,\infty)$ with density $f$, i.e. $F(x)=\int_0^x f(s)\d s$. Moreover, $m^{-1}:=\int_0^\infty xF(\d x)<\infty$.
\end{enumerate}

Let $\tau_1,\tau_2,\dots$ be a sequence of i.i.d. random variables whose cumulative distribution function is $F$. We consider the \emph{renewal process} that represents the arrival times of \emph{immigrants} and whose epochs are given by the \emph{renewals} $S_0=0$, $S_n=\tau_1+\dots+\tau_n$, $n\ge1$.
Each immigrant may have \emph{offspring}, and furthermore, each offspring individual can give place to offspring as well. To model this branching structure, let 
\begin{align}
\cbra{N_s^{\rbra{n}}(\cdot\mid t);\;t\in[0,\infty)}_{n\ge1},
\label{mblefamily}
\end{align} 
be a sequence of families of point processes which are considered to be i.i.d. and to be independent of $N_R$, such that $N_s^{(n)}(\cdot\mid t)$ has the same law as $N_s(\cdot\mid t)$, which for disjoint intervals $(a_1,b_1],$ $\dots,$ $(a_k,b_k]$,
\begin{align}
\label{PoissonComponent}
N_s((a_i,b_i]\mid t)\dist\text{Poi}\rbra{\int_{a_i}^{b_i}h(x-t)\d x}.
\end{align}
In other words, $N_s(\cdot\mid t)$ is an inhomogeneous Poisson process with intensity $h(x-t)\d x$, where we understand $h(x)=0$ for $x<0$. The matters of measurability of the indexing in \eqref{mblefamily} are treated in Section 3.2 of Hern\'andez--Yano \cite{IvanYanoConstruction2023}.

We now construct processes $N_c^{(n)}$, for $n\ge1$, from a superposition of the processes $N_s^{(n)}$ with the following recursive structure:
\begin{align}
\label{}
N_c^{(0)}(\cdot\mid t_0):=\delta_{t_0},\quad N_c^{\rbra{n+1}}(\cdot\mid t_0)=\sum_{t\in N_c^{\rbra{n}}(\cdot\mid t_0)}N_s^{\rbra{n+1}}(\cdot\mid t),
\end{align}
where $N_c^{(0)}(\cdot\mid t_0)$ is the original immigrant at $t_0$ and $N_c^{\rbra{n}}(\cdot\mid t_0)$ represents its $n$-th generation offspring. The total number of descendants in a cluster is given as,
\begin{align}
N_c(\cdot\mid t_0)=\sum_{n\ge 0} N_c^{(n)}(\cdot\mid t_0).
\end{align}
Hence, the total number of individuals is given by the superposition
\begin{align}
\label{RHPCluster}
N(\cdot)=\sum_{t_0\in N_R(\cdot)}N_c(\cdot\mid t_0)
\end{align}
In Hern\'andez--Yano \cite[Corollary 4.2]{IvanYanoConstruction2023}, it is shown that under assumptions \textbf{(A0)} and \textbf{(B0)}, $N(\cdot)$ correctly defines a point process $\cbra{T_i}_{i\ge0}$. Consider now, random variables
\begin{align}
D_i=\left\{\begin{array}{ll} 0\quad&\rbra{\text{if }T_i\in N_R(\cdot)},\\1&\rbra{\text{otherwise}}, \end{array}\right.
\end{align}
and the function $I(t)=\max\cbra{i;\;T_i\le t,D_i=0}$, $t\ge0$, which represents the index of the last immigrant up to time $t$. We construct the filtration $(\cF_t)_{t\ge0}$ by the augmentation of the natural filtration $(\cF^0_t)_{t\ge0}$ defined as
\begin{align}
\cF^0_t=\sigma\rbra{N_c^{(n)}\rbra{(a,b]}; \;0\le a\le b\le t, \;n=0,1,2,\dots},\quad t\ge0.
\end{align}
In this setting, using the uniqueness of intensities (up to a predictable modification)\cite[Sec II. T12]{BremaudQueues}, Hern\'andez--Yano \cite[Theorem 5.1]{IvanYanoConstruction2023} proved that the point process \eqref{RHPCluster} admits the $(\cF_t)$-intensity
\begin{align}
\lambda(t)=\mu\rbra{t-T_{I(t)}}+\int_0^{t} h(t-u)N(\d u).
\end{align}
Then, by setting $I$ and $(\cF_t)$ as above, we can define the RHP through its intensity.
\begin{Def}A point process $N$ is called a \emph{renewal Hawkes process} (RHP) if $N$ admits the $(\cF_t)$-intensity,
\begin{align}
\label{RHPInt}
\lambda(t)=\mu\rbra{t-T_{I(t)}}+\int_0^{t} h(t-u)N(\d u),
\end{align}
where $h$ satisfies \textbf{(A0)} and the \emph{hazard function} $\mu$ is a measurable function on $[0,\infty)$ defined as
\begin{align}
\label{DefinitionHazard}
\mu(t)=\frac{f(t)}{1-\int_0^t f(s)\d s}
\end{align}
for the probability density function $f$ in \textbf{(B0)}.
\end{Def}
Additionally, we introduce the following assumption:
\begin{enumerate}
\item[\textbf{(A1)}] The function $h$ is bounded and $h(t)\tend{}{t\rightarrow\infty}0$.
\end{enumerate}

We proceed with the statement of our main results. We have a law of large numbers for the RHP, in which we show that the mean number of arrivals can be consistently estimated as follows. 

\begin{Thm}
\label{TheoremLLN}
Assume \emph{\textbf{(A0, A1)}} and \emph{\textbf{(B0)}}. Then,
\begin{align}
\label{EquationLLN1}
\supre{}{v\in[0,1]}\absol{T^{-1}N(Tv)-v\frac{m}{1-\alpha}}\tend{a.s.}{T\rightarrow\infty}0.
\end{align}
\end{Thm}

The central limit theorem for the RHP takes the form,
\begin{Thm}
\label{TheoremCLT}
Under assumptions \emph{\textbf{(A0, A1)}} and \emph{\textbf{(B0)}}, if $\int_0^\infty x^2F(\d x)<\infty$, the convergence in distribution
\begin{align}
\rbra{\frac{1}{\sqrt{T}}\rbra{N(Tv)-\bE\sbra{N(Tv)}}}_{v\in[0,1]}\tend{d}{T\rightarrow\infty}\rbra{\sigma W(v)}_{v\in[0,1]},
\end{align}
holds in the Skorokhod topology, where $\;\rbra{W(v)}_{v\in[0,1]}$ is a standard Brownian motion and 
\begin{align}
\sigma^2=\frac{m\alpha}{(1-\alpha)^3}+\frac{m^3\Var\sbra{\tau}}{(1-\alpha)^2},
\end{align}
and $\tau$ is a random variable such that for $x\ge0$, $\bP(\tau\le x)=F(x)$.
\end{Thm}
And finally, we have the following result of asymptotic normality:
\begin{Cor}
\label{CorollaryAsymptoticNormality}
Under assumptions \emph{\textbf{(A0, A1)}} and \emph{\textbf{(B0)}}, if $\int_0^\infty x^qF(\d x)<\infty$ for some $q>\frac{5}{2}$ and $\int_0^\infty x^{r}h(x)\d x<\infty$ for $r>\frac{1}{2}$, the convergence in distribution
\begin{align}
\rbra{\frac{1}{\sqrt{T}}N(Tv)-v\frac{m}{1-\alpha}\sqrt{T}}_{v\in[0,1]}\tend{d}{T\rightarrow\infty}\rbra{\sigma W(v)}_{v\in[0,1]},
\label{AsympNorm}
\end{align}
holds in the Skorokhod topology, where $\sigma$ is the same as in Theorem \ref{TheoremCLT} and $(W(v))_{v\in[0,1]}$ is a standard Brownian motion. If additionally $F$ has an exponential moment, i.e., exists $\eta>0$ such that $\int_0^\infty e^{\eta x}F(x)\d x<\infty$, then, it is enough that $\int_0^\infty x^{r}h(x)\d x<\infty$, for some $r>0$ for the convergence in distribution \eqref{AsympNorm} to hold.
\end{Cor}
Notice that assumption \textbf{(A1)} is not required by Bacry--Delattre--Hoffmann--Muzy \cite{BacryLimit} for the derivation of Corollary \ref{CorBacry}. However, our integration condition on $h$ is more relaxed than the one assumed by them. In all of our results, the cases for the classical Hawkes process are recovered by taking $\tau$ as an exponential random variable of parameter $\mu$, or equivalently $F(x)=1-e^{-\mu x}$, $x\ge0$, which has an exponential moment, e.g. $\eta=\mu/2$.

We follow a martingale approach and make use of the renewal theorems to prove our results . The structure of the paper is the following. In Section \ref{SectionMean} we solve renewal-type equations to obtain expressions that are useful for the proofs of our results. Section \ref{SectionLLN} is dedicated to the proof of Theorem \ref{TheoremLLN}. Finally, in Section \ref{SectionCLT} we prove Theorem \ref{TheoremCLT} and Corollary \ref{CorollaryAsymptoticNormality}. In all instances it is indicated when the results for the classical Hawkes process can be retrieved as a particular case of the results in this work.

\section{Results for the mean number of arrivals}
\label{SectionMean}
In the following, the symbol of a measure $\nu$ on $[0,\infty)$ is used as well for its cumulative function $\nu(t)=\nu([0,t])$. Conversely, the symbol of a non-decreasing right-continuous function $\nu(t)$ on $[0,\infty)$ is used as well for its Stieltjes measure $\nu(\d t)$ such that $\nu(t)=\nu([0,t])$. Note that this abuse of symbols is standard (see, for example \cite[Sec. V]{Asmussen2003}). We also use the convention that the convolution between a function and a measure is a function, 
\begin{align}
&F^{*0}(x)=\delta_0(x),
\\
&F^{*(n+1)}(x)=F^{*n}*F(x)=\int_0^xF^{*n}(x-y)F(\d y),\quad x\ge0,\;n\ge0,
\\
&F^{*(n+1)}(x)=F*F^{*n}(x)=\int_0^xF(x-y)F^{*n}(\d y),\quad x\ge0,\;n\ge0.
\end{align}

The proofs for the results in this section involve the solution of \emph{renewal equations}. A renewal equation has the form
\begin{align}
\label{RenewalEquation}
Z(t)=z(t)+\int_0^t Z(t-u)F(\d u),\quad t\ge0,
\end{align}
where $z:[0,\infty)\rightarrow[0,\infty)$ is known, $Z:[0,\infty)\rightarrow[0,\infty)$ is unknown and $F$ is a given finite measure on $[0,\infty)$. The following Lemma gives conditions on the existence and uniqueness of the solution.
\begin{Lem}[\textbf{\cite[Lemma 4.1.I]{DVJ}}]
When $z$ is measurable and bounded on finite intervals, i.e. $\sup_{0\le t\le T}\absol{z(t)}<\infty$ for all $T<\infty$, the  renewal equation \eqref{RenewalEquation} has a unique
measurable solution that is also bounded on finite intervals, and it is given by
\begin{align}
\label{RenewalEquationSolution}
Z(t)=\Phi*z(t)=\int_0^t z(t-u)\Phi(\d u),\quad t\ge0,
\end{align}
where 
\begin{align}
\Phi(t):=\sum_{n\ge0}F^{*n}(t),\quad t\ge0.
\end{align}
\end{Lem}
In the case that the measure is absolutely continuous, we also have the following result.
\begin{Lem}[\textbf{\cite[Sec. V, Proposition 2.7]{Asmussen2003}}]
\label{RenewalDensity}
Let $F$ have a density with respect to the Lebesgue measure, i.e. $F(\d s)=f(s)\d s$, then
\begin{align}
\phi(t)=\sum_{n\ge1}f^{*n}(t),\quad t\ge0,
\end{align}
is a solution to the renewal equation
\begin{align}
\phi(t)=f(t)+\int_0^t f(t-u)F(\d u).
\end{align}
\end{Lem}

In our particular case, let $N$ be an RHP and consider its \emph{imbedded renewal process} 
\begin{align}
\cbra{S_0,S_1,\dots,}:=\cbra{T_i: D_i=0},
\end{align} 
to which corresponds the counting process
\begin{align}
N_R(t)=\sum_{i\ge0} 1_{\cbra{S_i\le t}},\quad t\ge0,
\end{align}
and the \emph{renewal function}
\begin{align}
\Phi(t):=\sum_{n\ge0}F^{*n}(t),\quad t\ge0.
\end{align}
It is well known in literature (c.f. \cite[Sec. V, Theorem 2.4]{Asmussen2003}) that if $F([0,\infty))\le 1$, then the renewal function is the expected number of renewals up to time $t$:
\begin{align}
\bE\sbra{N_R(t)}=\Phi(t).
\end{align}

We introduce then the following Lemma.
\begin{Lem}
\label{RHPLImbedded}
Let $N$ be an RHP whose imbedded renewal process $N_R$ has interarrival distribution $F$. Then,  the counting process $N_R(t)=\sum_{i\ge0} 1_{\cbra{S_i\le t}}$, $t\ge0$, admits the $(\cF_t)$-intensity $\mu(t-T_{I(t)})$. Moreover 
\begin{align}
\label{RHPL3E1}
\bE\sbra{\int_0^t\mu\rbra{s-T_{I(s)}}\d s}=\Phi(t).
\end{align}
\end{Lem}
\Proof{
From the cluster construction \eqref{RHPCluster}, it is clear that
\begin{align}
\label{EqualityIntensity}
\mu(t-T_{I(t)})=\mu(t-S_{N_R(t)-1}),\quad t\ge0,
\end{align}
which corresponds to the intensity of a renewal process with interarrival distribution $F$. Now, since the process $C(s)=1_{(0,t]}(s)$, $s\ge0$ is $(\cF_t)$-predictable, we can compute
\begin{align}
\bE\sbra{\int_0^t\mu(s-S_{N_R(s)-1})\d s}=\bE\sbra{\int_0^t N_R(\d s)}=\bE\sbra{N_R(t)}=\Phi(t).
\end{align}
Finally, from \eqref{EqualityIntensity}, the proof follows.
}
Let us note that if the imbedded renewal process has an interarrival distribution that satisfies \textbf{(B0)}, then the induced renewal measure $\Phi(\d t)$ can be decomposed (see for example Stone \cite{Stone1966} and Hern\'andez \cite{IvanMaximum2023}) as a sum of measures 
\begin{align}
\Phi=\Phi_1+\Phi_2,
\label{StonesDecomp}
\end{align} 
where $\Phi_2([0,\infty))<\infty$ and $\Phi_1$ is absolutely continuous with bounded density $\varphi_1$ that satisfies
\begin{align}
\varphi_1(t)\tend{}{t\rightarrow\infty}m.
\end{align}
In order to study the arrivals related to the self-exciting part of the process, we define the following function:
\begin{align}
\psi(t)=\sum_{n\ge1}h^{*n}(t),\quad t\ge0,
\end{align}
where $*$ denotes the convolution of functions in the usual sense:
\begin{align}
h^{*1}=h,\quad h^{*(n+1)}(t)=\int_0^t h^{*n}(t-s)h(s)\d s=\int_0^t h(t-s)h^{*n}(s)\d s.
\end{align}
We can now state the following Lemma.
\begin{Lem}
\label{LemmaExpN}
Assume \emph{\textbf{(A0)}} and \emph{\textbf{(B0)}}. For any $t\ge0$, the mean number of events $\bE\sbra{N(t)}$ is given as, 
\begin{align}
\label{EquationExpN1}
\bE\sbra{N(t)}=\Phi(t)+\int_0^t\psi(t-s)\Phi(s)\d s.
\end{align}
\end{Lem}
\Proof{
Let $t\ge0$. Since the process $C(s)=1_{(0,t]}(s)$, $s\ge0$ is $(\cF_t)$-predictable, from the property \eqref{PropertyIntensity} of the intensity and Lemma \ref{RHPLImbedded}, we have
\begin{align}
\bE\sbra{N(t)}=&\bE\sbra{\int_0^t\mu(s-T_{I(s)})\d s}+\bE\sbra{\int_0^t\int_0^s h(s-u)N(\d u)\d s}
\label{}\\
=&\Phi(t)+\bE\sbra{\int_0^t h(t-s)N(s)\d s}
\label{}\\
=&\Phi(t)+\int_0^t h(t-s)\bE\sbra{N(s)}\d s.
\end{align}
This is a renewal type integral equation for $\bE\sbra{N(t)}$. Since the renewal function is always finite, $\Phi$ is bounded on finite intervals, then the integral equation has a unique solution bounded on finite intervals given by
\begin{align}
\bE\sbra{N(t)}=\Phi(t)+\int_0^t \psi(t-s)\Phi(s)\d s,\quad t\ge0,
\end{align}
and this concludes the proof.
}

In the following, it will be useful to write the process
\begin{align}
\label{CenteredN}
X(t):=N(t)-\bE\sbra{N(t)},\quad t\ge0,
\end{align}
as a linear functional of the characteristic martingale $M(t)=N(t)-\int_0^t\lambda(s)\d s$, as it is done below.
\begin{Lem}
\label{LemmaLinearFunct}
Assume \emph{\textbf{(A0)}} and \emph{\textbf{(B0)}}. Set 
\begin{align}
A(t):=M(t)+\int_0^t \mu\rbra{s-T_{I(s)}}\d s-\Phi(t).
\end{align}
Then, for all $t\ge0$, the process $(X_t)_{t\ge0}$ satisfies,
\begin{align}
\label{EquationLinearFunct1}
X(t)=A(t)+\int_0^t\psi(t-s)A(s)\d s.
\end{align}
\end{Lem}
\Proof{
We have
\begin{align}
&X(t)=M(t)+\int_0^t\lambda(s)\d s-\bE\sbra{N(t)}.
\end{align}
We want to write $X(t)$ as a solution of a renewal equation involving $A(t)$. From the proof of Lemma \ref{LemmaExpN} we have
\begin{align}
&X(t)=M(t)+\int_0^t\lambda(s)\d s-\Phi(t)-\int_0^t h(t-s)\bE\sbra{N(s)}\d s
\label{}\\
=&M(t)+\int_0^t\mu(s-T_{I(s)})\d s+\int_0^t\int_0^sh(s-u)N(\d u)\d s-\Phi(t)-\int_0^th(t-s)\bE\sbra{N(s)}\d s
\label{}\\
=&M(t)+\int_0^t\mu(s-T_{I(s)})\d s+\int_0^th(t-s)N(s)\d s-\Phi(t)-\int_0^th(t-s)\bE\sbra{N(s)}\d s
\label{}\\
=&M(t)+\int_0^t\mu(s-T_{I(s)})\d s-\Phi(t)+\int_0^th(t-s)X(s)\d s.
\label{}\\
=&A(t)+\int_0^th(t-s)X(s)\d s.
\end{align}
Since the function $\mu$ is locally integrable, then the process $(A_t)_{t\ge0}$ is a.s. bounded on finite intervals, therefore we have a solution for $X(t)$ given by
\begin{align}
X(t)=A(t)+\int_0^t\psi(t-s)A(s)\d s,\quad t\ge0.
\end{align}
The proof is complete.
}

\section{Law of large numbers}
\label{SectionLLN}

In preparation for the proof of Theorem \ref{TheoremLLN}, we need the following Lemmata.

\begin{Lem}
\label{LemmaLLNSpeedConvergence}
Assume \emph{\textbf{(A0, A1)}} and \emph{\textbf{(B0)}}. Then,
\begin{align}
\label{RHPLLLN1E0}
\supre{}{v\in[0,1]}\absol{T^{-1}\bE\sbra{N(Tv)}-v\frac{m}{1-\alpha}}\tend{}{T\rightarrow\infty}0.
\end{align}
If additionally we assume $\int_0^\infty x^{r}h(x)\d x<\infty$ for some $r>0$ and $\int_0^\infty x^{q}F(\d x)<\infty$ for some $q>2$, we have for any $0\le p<\min\{1,r,q-2\}$ that
\begin{align}
\label{RHPLLLN1E1}
T^p\supre{}{v\in[0,1]}\absol{T^{-1}\bE\sbra{N(Tv)}-v\frac{m}{1-\alpha}}\tend{}{T\rightarrow\infty}0.
\end{align}
If instead $F$ has an exponential moment, i.e., exists $\eta>0$ such that $\int_0^\infty e^{\eta x}F(x)\d x<\infty$, and $\int_0^\infty x^{r}h(x)\d x<\infty$, for some $r>0$, then the convergence \eqref{RHPLLLN1E1} holds for all $0<p<1$.
\end{Lem}
\Proof{Let us recall Lemma \ref{LemmaExpN},
\begin{align}
\bE\sbra{N(t)}=\Phi(t)+\int_0^t\psi(t-s)\Phi(s)\d s,
\end{align}
 Stone's decomposition \eqref{StonesDecomp}, and the fact that $v\frac{m}{1-\alpha}=v\rbra{m+\frac{m\alpha}{1-\alpha}}$, we have
\begin{align}
&-T^p\rbra{T^{-1}\bE\sbra{N(Tv)}-v\frac{m}{1-\alpha}}
\label{}\\
=&T^p\sbra{\rbra{vm-\frac{\Phi(Tv)}{T}}+\rbra{v\frac{m\alpha}{1-\alpha}-\frac{\int_0^{Tv}\psi(Tv-s)\Phi(s)\d s}{T}}}.
\label{RHPLLLN1E2}
\end{align}
We define the function
\begin{align}
G(t):=\int_0^t\psi(t-s)\Phi(s)\d s,\quad t\ge0.
\end{align}
By changing the order of integration (note that $\Phi(t)=\int_0^t\Phi(\d s)$) we can rewrite $G$ as
\begin{align}
G(t)=\int_0^t\Psi(t-s)\Phi(\d s),\quad t\ge0,
\end{align}
where
\begin{align}
H(t):=\int_0^th(s)\d s,\quad \Psi(t):=\sum_{n\ge0}H^{*n}(t),\quad t\ge0.
\end{align}
Notice that $\Psi$ can be defined for all $t\ge0$ because it can be seen as the renewal function for the defective distribution $H$ with density $h$ (see for example Lemma \ref{RenewalDensity} above). It then holds that
\begin{align}
\Psi(t)=\int_0^t\psi(s)\d s,\quad t\ge0.
\end{align}
From the relation $\rbra{f*g}^\prime=f^\prime*g$, we can deduce
\begin{align}
G^\prime(t)=\int_0^t\psi(t-s)\Phi(\d s),\quad t\ge0,
\end{align}
in the sense of the Radon--Nikodym derivative, i.e.
\begin{align}
G(t)=\int_0^t G^\prime(s)\d s.
\end{align}
Let us analyze the asymptotics of $\psi$. For this, note that we can write $\psi$ as the solution to the renewal equation, 
\begin{align}
\psi=h+\psi*H.
\label{RHPLLLNERenforh}
\end{align}
Since $H(\infty):=\lim_{t\rightarrow\infty}H(t)=\alpha<1$, equation \eqref{RHPLLLNERenforh} is a defective renewal equation \cite[Sec V. Eq.(2.1)]{Asmussen2003} with solution $\psi=\sum_{n\ge1}h^{*n}$. In the case of a defective renewal equation, it holds that (c.f. V.7.4 in \cite{Asmussen2003}),
\begin{align}
\psi(t)\tend{}{t\rightarrow\infty}\frac{h(\infty)}{1-H(\infty)}=\frac{0}{1-\alpha}=0.
\label{EquationConvergencePsi}
\end{align}
Write $\norm{h}_\infty:=\supre{}{t\ge0}\absol{h(t)}$. Boundedness of $\psi$ can be seen from
\begin{align}
&h^{*n}*h(t)=\int_0^t h^{*n}(t-s)h(s)\d s\le\alpha\norm{h^{*n}}_\infty,
\label{}\\
&\norm{h^{*n}}_\infty\le\alpha^{n-1}\norm{h}_\infty,
\label{}\\
&\norm{\psi}_\infty=\norm{\sum_{n\ge1}h^{*n}}_\infty\le \frac{\norm{h}_\infty}{1-\alpha}.
\label{EquationPsiBounded}
\end{align}
Then from \eqref{EquationConvergencePsi}, \eqref{EquationPsiBounded}, and the Key Renewal Theorem in the case of an absolutely continuous interarrival distribution (c.f. Corollary VII.1.3 in \cite{Asmussen2003}), 
\begin{align}
&G^\prime(\infty):=\lim_{t\rightarrow\infty}G^\prime(t)= m\int_0^\infty\psi(s)\d s
\label{}\\
=&m\sum_{n\ge1}\int_0^\infty h^{*n}(s)\d s=m\sum_{n\ge1}\alpha^n=\frac{m\alpha}{1-\alpha}.
\label{}
\end{align}
The reasoning above allows us to rewrite \eqref{RHPLLLN1E2} as
\begin{align}
&T^p\sbra{\rbra{vm-\frac{\Phi_2(Tv)+\int_0^{Tv}\varphi_1(s)\d s}{T}}+\rbra{vG^\prime(\infty)-\frac{\int_0^{Tv}G^\prime(s)\d s}{T}}}
\label{}\\
=&\sbra{\rbra{\frac{\int_0^{Tv}m-\varphi_1(s)\d s}{T^{1-p}}}+\rbra{\frac{\int_0^{Tv}G^\prime(\infty)-G^\prime(s)\d s}{T^{1-p}}}-\frac{\Phi_2(Tv)}{T^{1-p}}}.
\end{align}
Now we notice that
\begin{align}
&T^p\supre{}{v\in[0,1]}\absol{\rbra{T^{-1}\bE\sbra{N(Tv)}-v\frac{m}{1-\alpha}}}
\label{}\\
\le &\supre{}{v\in[0,1]}\absol{\frac{\int_0^{Tv}m-\varphi_1(s)\d s}{T^{1-p}}} +\supre{}{v\in[0,1]}\absol{\frac{\int_0^{Tv}G^\prime(\infty)-G^\prime(s)\d s}{T^{1-p}}}+\frac{\Phi_2([0,\infty))}{T^{1-p}}
\label{vanishinglln}\\
\le&\frac{\int_0^{T}\absol{m-\varphi_1(s)}\d s}{T^{1-p}} +\frac{\int_0^{T}\absol{G^\prime(\infty)-G^\prime(s)}\d s}{T^{1-p}}+\frac{\Phi_2([0,\infty))}{T^{1-p}}.
\end{align}
If we assume \textbf{(B0)} and take $p=0$, the result \eqref{RHPLLLN1E0} follows immediately from the finiteness of the measure $\Phi_2$ and the convergence of $G^\prime(t)\tend{}{t\rightarrow\infty}G^\prime(\infty)$ and $\varphi_1(t)\tend{}{t\rightarrow\infty}m$.

Let us now assume that $\int_0^\infty x^{r}h(x)\d x<\infty$ for some $r>0$ and $\int_0^\infty x^{q}F(\d x)<\infty$ for some $q>2$, and let us study the case for $0<p<\min\{1,q-2\}$. We want to study the rates of convergence of $\varphi_1$ and $G^\prime$. Hern\'andez \cite[Lemma 3.2]{IvanMaximum2023} showed that given $\int_0^\infty x^q F(\d x)<\infty$, it holds that 
\begin{align}
\varphi_1(s)=m+o(s^{1-q})\quad \text{as } s\rightarrow\infty. 
\end{align}
Hence,
\begin{align}
&\frac{\int_0^{T}\absol{m-\varphi_1(s)}\d s}{T^{1-p}}=\frac{O(T^{2-q})}{T^{1-p}}\tend{}{T\rightarrow\infty}0.
\label{}
\end{align}
To study the integral with $G^\prime$, we first recall that for any $n\ge1$,
\begin{align}
\absol{\sum_{k=1}^{n}x_k}^r\le\left\{\begin{array}{rc} n^{r-1}\sum_{k=1}^{n}\absol{x_k}^r,& r>1,\\\sum_{k=1}^{n}\absol{x_k}^r,&r\le1, \end{array}\right.
\end{align}
from which it follows that for $r>0$ 
\begin{align}
\int_0^\infty x^{r}\psi(x)\d x=&\sum_{n=1}^\infty \int_0^\infty x^{r} h^{*n}(x)\d x
\label{}\\
=&\sum_{n=1}^\infty\int_0^\infty\cdots\int_0^\infty (x_1+\dots+x_n)^{r}h(x_1)\cdots h(x_n)\d x_1\cdots\d x_n
\label{}\\
\le&\left\{\begin{array}{rc} \sum_{n=1}^\infty n^{r}\alpha^{n-1}\int_0^\infty x^r h(x)\d x,& r>1,\\ \sum_{n=1}^\infty n\alpha^{n-1}\int_0^\infty x^r h(x)\d x,&r\le1, \end{array}\right.<\infty.
\end{align}
We will now show that  
\begin{align}
\frac{1}{T^{1-p}}\int_1^{T}\absol{G^\prime(\infty)-G^\prime(s)}\d s\tend{}{T\rightarrow\infty}0.
\label{EquationFinitenessG}
\end{align}
We can proceed along the lines of Hern\'andez \cite[Proof of Theorem 1.5]{IvanMaximum2023} from which we have the bound
\begin{align}
\label{boundlln}
\absol{G^\prime(\infty)-G^\prime(x)}\le &C\int_{x/2}^\infty\psi(y)\d y+\int_{x/2}^{x}\absol{\varphi_1(y)-m}\d y+\int_0^{x}\psi(x-y)\Phi_2(\d y),
\end{align}
and the constant $C$ can be taken as $m+\norm{\varphi_1-m}_{\infty}<\infty$. Looking at the first integral on the RHS of \eqref{boundlln}, we notice that
\begin{align}
\int_1^{T}\rbra{\int_{x/2}^\infty \psi(y)\d y}\d x\le&\int_1^{T}\rbra{\int_{x/2}^\infty \frac{(2y)^{r}}{x^{r}}\psi(y)\d y}\d x
\label{}\\
=&\int_1^{T}x^{-r}\d x\int_{x/2}^\infty (2y)^{r}\psi(y)\d y
\label{}\\
\le&\int_1^{T}x^{-r}\d x\int_{0}^\infty (2y)^{r}\psi(y)\d y=O(T^{1-r})\quad\text{as }T\rightarrow\infty.
\end{align}
For the second summand of \eqref{boundlln}, we obtain the following
\begin{align}
\int_1^T\rbra{\int_{x/2}^{x}\absol{\varphi_1(y)-m}\d y }\d x\le&\int_{1/2}^{T}\rbra{\int_{y}^{2y}\absol{\varphi_1(y)-m}\d x} \d y
\label{}\\
=&\int_{1/2}^{T}y\absol{\varphi_1(y)-m} \d y
\label{}\\
=&O(T^{3-q})\quad\text{as }T\rightarrow\infty,
\end{align}
which follows once again from $\varphi_1(y)=m+o(y^{1-q})$ as $y\rightarrow\infty$. Finally, using Fubini's theorem we see that
\begin{align}
\int_0^{\infty}\int_0^x \psi(x-y)\Phi_2(\d y)\d x=\int_0^\infty\psi(u)\d u\int_0^\infty\Phi_2(\d y)<\infty.
\end{align}
Therefore, for the values of $p$ considered, all terms in \eqref{vanishinglln} vanish as $T\rightarrow\infty$. 
Finally, if we assume that it exists $\eta>0$ such that $\int_0^\infty e^{\eta x}F(x)\d x<\infty$, then it holds (c.f. \cite[Sec. VII, Theorem 2.10]{Asmussen2003}) that there exists $\epsilon>0$ such that
\begin{align}
\varphi_1(s) = m+o(e^{-\epsilon t})\quad\text{as }s\rightarrow\infty
\end{align}
Using the exponential decay in the above expressions concludes the proof.
}
\textbf{Remark.} The result about exponential decay in \cite[Sec. VII, Theorem 2.10]{Asmussen2003} is stated with $O(e^{-\epsilon t})$, however, following the justification on \cite[Sec. VII, page 191]{Asmussen2003}, we see that in fact this decay is with $o(e^{-\epsilon t})$ asymptotics.

The next result treats the asymptotic behavior of the renewal process part.
\begin{Lem}
\label{RHPLLNIntTerm}
Under \emph{\textbf{(B0)}}, we have almost surely that
\begin{align}
T^{-1}\supre{}{t\le T}\absol{\int_0^t\mu(s-T_{I(s)})\d s-\Phi(t)}\tend{}{T\rightarrow\infty}0.
\end{align}
\end{Lem}
\Proof{
We rewrite the integral term as
\begin{align}
\int_0^t\mu(s-T_{I(s)})\d s=&\int_0^t\mu(s-S_{N_R(s)-1})\d s
\label{}\\
=&\sum_{j=1}^{N_R(t)-1}\int_{S_{{j-1}}}^{S_{j}}\mu(s-S_{{j-1}})\d s+\int_{S_{{N_R(t)-1}}}^{t}\mu(s-S_{{N_R(t)-1}})\d s
\label{}\\
=&\sum_{j=1}^{N_R(t)-1}\xi_j+\int_{S_{{N_R(t)-1}}}^{t}\mu(s-S_{{N_R(t)-1}})\d s,
\end{align}
where the $\xi_j$ are i.i.d. random variables. To compute their mean we use the definition of intensity of a point process \eqref{DefIntensity} and the fact that the process $1_{(S_{{j-1}},S_{j}]}(t)$ is predictable, and we obtain
\begin{align}
\bE\sbra{\xi_j}=&\bE\sbra{\int_{S_{{j-1}}}^{S_{j}}\mu(s-S_{{j-1}})\d s}
\label{}\\
=&\bE\sbra{\int_{0}^{\infty}1_{(S_{{j-1}},S_{j}]}(s)\mu(s-S_{{j-1}})\d s}
\label{}\\
=&\bE\sbra{\int_{0}^{\infty}1_{(S_{{j-1}},S_{j}]}(s)N_R(\d s)}
\label{}\\
=&\bE\sbra{N_R\rbra{{(S_{{j-1}},S_{j}]}}}=1,
\end{align}
for all $1\le j\le N_R(t)-1$. On the one hand, from the Law of Large Numbers for $\xi_j$ we have
\begin{align}
\frac{1}{n}\sum_{j=1}^n\xi_j\tend{a.s.}{}1,
\end{align}
while from the LLN for the inter-arrival times of $N_R(\cdot)$ we have
\begin{align}
\frac{1}{n}S_{n}=\frac{1}{n}\sum_{j=1}^n\tau_j\tend{a.s.}{}\frac{1}{m}.
\end{align}
Combining these two facts we obtain that
\begin{align}
&\frac{1}{S_{n}}\int_{0}^{S_{n}}\mu(s-T_{I(s)})\d s=\frac{1}{S_{n}}\sum_{j=1}^{n}\int_{S_{{j-1}}}^{S_{j}}\mu(s-S_{{j-1}})\d s=\frac{1}{S_{n}}\sum_{j=1}^n\xi_j\tend{a.s.}{}m
\end{align}
Since for $t>0$ we may take $n=n(t)$ such that $S_{{n-1}}<t\le S_{{n}}$,
\begin{align}
R_t:=\frac{1}{t}\int_{0}^{t}\mu(s-T_{I(s)})\d s\tend{a.s.}{t\rightarrow\infty}m,
\end{align}
because,
\begin{align}
\frac{1}{S_{{n-1}}+\tau_{n}}\int_{0}^{S_{{n-1}}}\mu(s-T_{I(s)})\d s\le R_t\le \frac{1}{S_{{n}}-\tau_{n-1}}\int_{0}^{S_{n}}\mu(s-T_{I(s)})\d s,
\label{}
\end{align}
and the limit on each side is equal to $m$ a.s.

Furthermore, from the elementary renewal theorem, i.e.
\begin{align}
\frac{\Phi(t)}{t}\tend{}{t\rightarrow\infty}m,
\end{align}
and the fact that
\begin{align}
&T^{-1}\supre{}{t\le T}\absol{\int_0^t\mu(s-T_{I(s)})\d s-\Phi(t)}
\\\label{}
\le &T^{-1}\supre{}{t\le T}\absol{\int_0^t\mu(s-T_{I(s)})\d s-m}+T^{-1}\supre{}{t\le T}\absol{m-\Phi(t)},
\end{align}

we obtain the desired conclusion.
}

We can proceed with the proof of Theorem \ref{TheoremLLN} that is very similar to the proof wrote by Bacry--Delattre--Hoffmann--Muzy \cite[Theorem 1]{BacryLimit}, but which we include with slightly more detail for the sake of completion of this paper.
\Proof[Proof of Theorem \ref{TheoremLLN}]{
We use \eqref{RHPLLLN1E0} of Lemma \ref{LemmaLLNSpeedConvergence}, then it suffices to prove that
\begin{align}
T^{-1}\supre{}{v\in[0,1]}\absol{N(Tv)-\bE\sbra{N(Tv)}}\tend{a.s.}{T\rightarrow\infty}0.
\end{align}
From Lemma \ref{LemmaLinearFunct} we know that for $T\ge0$ and $v\in[0,1]$
\begin{align}
X(Tv)=A(Tv)+\int_0^{Tv}\psi(Tv-s)A(s)\d s,
\end{align}
with $A(Tv)=M(Tv)+\int_0^{Tv} \mu\rbra{s-T_{I(s)}}\d s-\Phi(Tv)$. Then
\begin{align}
\supre{}{v\in[0,1]}\absol{X(Tv)}\le&\supre{}{t\le T}\absol{A(t)}+\supre{}{t\le T}\int_0^t\psi(t-s)\absol{A(s)}\d s
\label{}\\
\le&\supre{}{t\le T}\absol{A(t)}\rbra{1+\int_0^\infty\psi(s)\d s},
\end{align} 
where we note that $\psi(\cdot)$ is integrable. Now, shifting attention to $\supre{}{t\le T}\absol{A(t)}$, we obtain the bound
\begin{align}
\supre{}{t\le T}\absol{A(t)}\le& \supre{}{t\le T}\absol{M(t)}+\supre{}{t\le T}\absol{\int_0^t\mu(s-T_{I(s)})\d s-\Phi(t)}.
\end{align}
Consider the characteristic martingale $M(t)$, and define the martingale
\begin{align}
Z(t)=\int_{(0,t]}\frac{1}{s+1}\d M_s.
\end{align}
We compute its quadratic variation
\begin{align}
\sbra{Z,Z}_t=&\sum_{0<u\le t}\rbra{Z(u)-Z(u-)}^2
\label{}\\
=&\sum_{0<u\le t}\rbra{\int_{(0,u]}\frac{1}{s+1}\d M_s-\int_{(0,u)}\frac{1}{s+1}\d M_s}^2
\label{}\\
=&\sum_{0<u\le t}\rbra{\frac{1}{u+1}\rbra{N(u)-N(u-)}}^2
\label{}\\
=&\int_{(0,t]}\frac{1}{\rbra{s+1}^2}N(\d s).
\end{align}
Now, by integration by parts, we obtain
\begin{align}
\int_0^t\frac{1}{\rbra{u+1}^2}N(\d u)-\frac{N(t)}{\rbra{t+1}^2}=&\int_0^t\sbra{\frac{1}{\rbra{u+1}^2}-\frac{1}{\rbra{t+1}^2}}N(\d u)
\label{}\\
=&\int_0^t\int_u^t\frac{2}{\rbra{s+1}^3}\d sN(\d u)
\label{}\\
=&2\int_0^t\int_0^s\frac{N(\d u)}{\rbra{s+1}^3}\d s
\label{}\\
=&2\int_0^t\frac{N(s)}{\rbra{s+1}^3}\d s.
\end{align}
From the previous equalities, we have
\begin{align}
&\bE\sbra{\int_0^t\frac{1}{\rbra{s+1}^2}N(\d s)}=2\bE\sbra{\int_0^t\frac{N(s)}{\rbra{s+1}^3}\d s}+\frac{\bE\sbra{N(t)}}{(t+1)^2}.
\end{align}
We can analyze the second term on the RHS by using Lemma \ref{LemmaExpN} and the increasingness of $\Phi(\cdot)$,
\begin{align}
\frac{\bE\sbra{N(t)}}{(t+1)^2}=\frac{\Phi(t)}{\rbra{t+1}^2}+\int_0^t\frac{\psi(t-s)}{\rbra{t+1}^2}\Phi(s)\d s\le \frac{\Phi(t)}{\rbra{t+1}^2}\cbra{1+\int_0^t\psi(t-s)\d s}\tend{}{t\rightarrow\infty}0.
\end{align}
Using the monotone convergence theorem and the elementary renewal theorem, we obtain, 
\begin{align}
&\bE\sbra{\int_0^\infty\frac{1}{\rbra{s+1}^2}N(\d s)}=2\bE\sbra{\int_0^\infty\frac{N(s)}{\rbra{s+1}^3}\d s}
\label{}\\
=&2\int_0^\infty\frac{\bE\sbra{N(s)}}{\rbra{s+1}^3}\d s=2\int_0^\infty\frac{\Phi(s)}{\rbra{s+1}^3}\d s+2\int_0^\infty\frac{\int_0^s\psi(s-u)\Phi(u)\d u}{\rbra{s+1}^3}\d s
\label{}\\
\le&2\int_0^\infty\frac{\Phi(s)}{\rbra{s+1}^3}\d s+2\int_0^\infty\frac{\Phi(s)}{\rbra{s+1}^3}\d s\int_0^\infty\psi(u)\d u<+\infty.
\end{align}
This tells us that $Z(\cdot)$ is a martingale bounded in $L^2(\bP)$, therefore, by the martingale convergence theorem $\lim_{t\rightarrow\infty}Z(t)$ exists and is finite a.s. Let us recall that $M(0)=0$, and consider
\begin{align}
\int_0^t Z(s)\d s=&\int_0^t\int_0^s\frac{\d M_u}{u+1}\d s
\label{}\\
=&\int_0^t\frac{t-u}{u+1}\d M_u
\label{}\\
=&(t+1)\int_0^t\frac{\d M_u}{u+1}-\int_0^t\frac{u+1}{u+1}\d M_u
\label{}\\
=&(t+1)Z(t)-M(t).
\end{align}
Furthermore, from the finiteness of the limit of $Z(t)$ it holds that
\begin{align}
\frac{1}{t+1}M(t)=Z(t)-\frac{1}{t+1}\int_0^tZ(s)\d s\tend{a.s.}{t\rightarrow\infty}0.
\end{align}
Finally, we show that the convergence is uniform in $v\in[0,1]$. Let $0<\epsilon<1$. For $0\le v<\epsilon$, we have,
\begin{align}
&\absol{\frac{M(Tv)}{T}}=\absol{\frac{M(Tv)}{Tv}v}\le\epsilon\supre{}{0<t<\infty}\absol{\frac{M(t)}{t}}.
\end{align}
Meanwhile, for $\epsilon\le v\le1$, 
\begin{align}
&\absol{\frac{M(Tv)}{T}}\le\supre{}{T\epsilon<t<\infty}\absol{\frac{M(t)}{t}}.
\end{align}
Hence, we obtain,
\begin{align}
&\limsup_{T\rightarrow\infty}\rbra{\supre{}{v\in[0,1]}\absol{\frac{M(Tv)}{T}}}\le\supre{}{0<t<\infty}\absol{\frac{M(t)}{t}}\epsilon.
\end{align}
Since this was done for an arbitrary $\epsilon$, we conclude the desired uniform convergence.
}
It should be pointed out that these results depend on the asymptotic and integration properties of the inter arrival distribution for the renewals and the excitation function $h$. A possible generalization to a multivariate RHP could for example treat the case where all the immigration processes are independent from each other, their distributions have the required finite moment properties for each component, and the excitation functions between components satisfy the asymptotic and integration conditions stated above.

\section{Central limit theorem}
\label{SectionCLT}
In this section, we make use of the cluster construction of the RHP to derive a central limit theorem.
\begin{Lem}
\label{RHPLCLTSumprocess}Assume \emph{\textbf{(A0, A1)}}, \emph{\textbf{(B0)}} and $\int_0^\infty x^2F(\d x)<\infty$. For $t\ge0$, set 
\begin{align}
X^{(T)}(v):=\frac{N(vT)-\bE\sbra{N(vT)}}{\sqrt{T}},\quad T > 0
\end{align}
Then, 
\begin{align}
X^{(T)}(v)\tend{d}{T\rightarrow\infty}N(0,\sigma^2v),
\end{align}
where
\begin{align}
\sigma^2=\frac{m\alpha}{(1-\alpha)^3}+\frac{m^3\Var\sbra{\tau}}{(1-\alpha)^2}.
\end{align}
Moreover, for $0<u<v<1$,
\begin{align}
\Cov\sbra{X^{(T)}(u),X^{(T)}(v)}\tend{}{T\rightarrow\infty}\sigma^2 u.
\end{align}
\end{Lem}
\Proof{
Let $N$ be an RHP and we consider its cluster representation \eqref{RHPCluster}. We denote as 
\begin{align}
W=N_c\rbra{[0,\infty)\mid 0},
\end{align}
which represents the total size of a Galton--Watson process with offspring distribution Poisson($\alpha$). This random variable has mean $\bE\sbra{W}=\frac{1}{1-\alpha}$ and variance $\Var\sbra{W}=\frac{\alpha}{(1-\alpha)^3}$ (c.f. \cite[Section 1.3]{Harris1963}). Consider the sequence $\cbra{W_i}_{i]\ge0}$ of i.i.d random variables with the same distribution as $W$.  We claim that for large $T$, the points coming from clusters close to $T$, but that lie outside $[0,T]$ become negligible in comparison to $N(T)$. Indeed, let us define
\begin{align}
Y_u:=N_c((u,\infty)), \quad u\ge0,
\end{align}
which represents the descendants that arrive after a delay $u$. We can then represent the \emph{escaping points} that fall outside the interval $[0,T]$ for any $T\ge0$ as
\begin{align}
N_A([0,T])=\sum_{k=1}^{N_R(T)}Y_{T-S_k}.
\end{align}
We now use a result of Reynaud-Bouret--Roy\cite[Proposition 1.3]{ReynaudBouret} to find constants $C,\beta$ and $z_0\in [0,\alpha-\log\alpha-1]$ such that 
\begin{align}
\bE\sbra{e^{z_0 Y_u}} \le \exp\rbra{Ce^{-\beta u}},\quad \bE\sbra{Y_u} \le C_1 e^{-\beta u}, \quad C_1:=\tfrac{Cz_0}{\beta}.
\label{RBineq}
\end{align}
Now, we will choose a slowly growing window of observation before $T$. For this set
\begin{align}
R_T \;:=\; \Bigl\lceil \tfrac{2}{\beta}\log T \Bigr\rceil,\qquad T\ge 2,
\label{RTdef}
\end{align}
so that $R_T<T$, $R_T\to\infty$ and $R_T/T\tend{}{T\rightarrow\infty} 0$. We can then split the escaping points into those originated by immigrants before the observation window and those originated within it:
\begin{align}
N_A[0,T] \;=\; A_T + B_T,
\end{align}
where
\begin{align}
A_T := \sum_{S_k\le T-R_T} Y_{T-S_k}, \qquad B_T := \sum_{T-R_T< S_k\le T} Y_{T-S_k}.
\end{align}
Using \eqref{RBineq} and \eqref{RTdef},
\begin{align}
\bE\sbra{A_T}\le C_1 e^{-\beta R_T}\, N_R(T-R_T)\le \frac{C_2}{T^{2}},
  \qquad C_2 := m C_1.
\end{align}
Hence, by Markov's inequality and Borel–Cantelli,
\begin{equation}
\frac{A_T}{T} \tend{}{T\rightarrow\infty} 0\quad {\text{a.s.}}
\label{ATvanish}
\end{equation}
For the second term, let $M_T:=N_R(T)-N_R(T-R_T)$ be the number of ancestors in $(T-R_T,T]$. By standard arguments of renewal theory, we get
\begin{align}
\frac{M_T}{R_T}\;\tend{}{T\rightarrow\infty}m\quad {\text{a.s.}},
\end{align}
so almost surely $M_T\le 2 m R_T$ for all $T$ large enough. 

For any $u\le R_T$ the bound \eqref{RBineq} gives
\begin{align}
\bE\sbra{e^{z_0 Y_u}} \;\le\; e^{C},
\end{align}
hence
\begin{align}
\bE\sbra{e^{z_0 B_T}}=  \bE\sbra{\rbra{\bE[e^{z_0 Y_u}]}^{M_T}}\le  e^{C M_T}\le\exp\rbra{C'R_T} =  T^{C'2\beta}.
\end{align}
Chernoff’s bound with parameter $z_0$ yields
\begin{align}
\bP\rbra{B_T\ge\varepsilon T}\le \exp\rbra{-z_0\varepsilon T + C''\log T},
\end{align}
whose series over $T$ converges, so once again, using Borel–Cantelli gives
\begin{align}
\frac{B_T}{T} \tend{}{T\rightarrow\infty} 0\quad {\text{a.s.}}.
\label{BTvanish}
\end{align}
Then, by combining the almost-sure limits \eqref{ATvanish} and \eqref{BTvanish}, we get
\begin{align}
\frac{N_A[0,T]}{T}=\frac{A_T}{T}+\frac{B_T}{T}\tend{}{T\rightarrow\infty} 0\quad {\text{a.s.}}
\end{align}

From the previous argument we learned that we can consider all the clusters $W_i$ that were originated by immigrants up to time $T$. The escaping points will be negligible after dividing by $T$ and taking the limit. Thus, we can write
\begin{align}
X^{(T)}(v)=&\frac{\sum_{i=0}^{N_R(vT)}W_i-\bE\sbra{\sum_{i=0}^{N_R(vT)}W_i}}{\sqrt{T}}
\label{}\\
=&\frac{\sum_{i=0}^{N_R(vT)}\rbra{W_i-\bE\sbra{W_i}}}{\sqrt{T}}+\frac{\bE\sbra{W}\rbra{N_R(vT)-\bE\sbra{N_R(vT)}}}{\sqrt{T}}.
\end{align}
Let us analyze each term separately. First we rewrite the first term as
\begin{align}
\frac{\sum_{i=0}^{N_R(vT)}\rbra{W_i-\bE\sbra{W_i}}}{\sqrt{N_R(vT)}}\frac{\sqrt{N_R(vT)}}{\sqrt{T}}
\end{align}
This is a random sum of i.i.d random variables of variance $\frac{\alpha}{(1-\alpha)^3}$. We use the elementary renewal theorem to find $\sqrt{N_R(vT)}/\sqrt{T}\tend{\text{a.s.}}{n\rightarrow\infty} \sqrt{mv}$ and Anscombe's theorem (c.f. \cite[p.216]{Chung1974}) to find that
\begin{align}
\frac{\sum_{i=0}^{N_R(vT)}\rbra{W_i-\bE\sbra{W_i}}}{\sqrt{T}\sqrt{N_R(vT)}}\frac{\sqrt{N_R(vT)}}{\sqrt{T}}\tend{d}{T\rightarrow\infty}N\rbra{0,\frac{mv\alpha}{(1-\alpha)^3}}.
\end{align}
Notice now that the above term has the same limit in distribution as
\begin{align}
\frac{\sum_{i=0}^{\bE\sbra{N_R(vT)}}\rbra{W_i-\bE\sbra{W_i}}}{\sqrt{T}}.
\label{CLTProofAux1}
\end{align}

Indeed, let us consider the expectation,
\begin{align}
&\bE\sbra{\rbra{\frac{1}{\sqrt{T}}\sum_{i=0}^{N_R(vT)}\rbra{W_i-\bE\sbra{W_i}}-\frac{1}{\sqrt{T}}\sum_{i=0}^{\bE\sbra{N_R(vT)}}\rbra{W_i-\bE\sbra{W_i}}}^2}
\label{}\\
=&\frac{1}{T}\bE\sbra{\absol{N_R(vT)-\bE\sbra{N_R(vT)}}\Var\sbra{W}}
\label{}\\
\le&\Var\sbra{W}\rbra{\bE\sbra{\absol{\frac{N_R(vT)}{T}-mv}}+\absol{\frac{\bE\sbra{N_R(vT)}}{T}-mv}}\tend{}{T\rightarrow\infty}0.
\end{align}

For the second term, we can directly apply the central limit theorem for the renewal process \cite[Sec V. Proposition 6.3]{Asmussen2003} to see that
\begin{align}
\frac{\bE\sbra{W}\rbra{N_R(vT)-\bE\sbra{N_R(vT)}}}{\sqrt{T}}\tend{d}{T\rightarrow\infty}N\rbra{0,\frac{m^3\Var\sbra{\tau}v}{(1-\alpha)^2}}.
\label{CLTProofAux2}
\end{align}
Since \eqref{CLTProofAux1} is independent from \eqref{CLTProofAux2}, then the limit is the sum of the limiting normal random variables, i.e., $X^{(T)}(v)\tend{d}{T\rightarrow\infty}N(0,\sigma^2)$. Now, let $0<u<v<1$. We begin by noticing that
\begin{align}
&X^{(T)}(v) = X^{(T)}(u) + \sum_{i=N_R(uT)}^{N_R(vT)} W_i
\label{}\\
&\Cov\sbra{X^{(T)}(u),X^{(T)}(v)}=\Cov\sbra{X^{(T)}(u), X^{(T)}(u) + \sum_{i=N_R(uT)}^{N_R(vT)} W_i }.
\end{align}
If we expand the covariance and use the independence of the clusters that arrived before $u$ from the clusters that arrived  in the interval $[u,v)$, we get that
\begin{align}
\Cov\sbra{X^{(T)}(u), X^{(T)}(v)} =& \Cov\sbra{X^{(T)}(u)), X^{(T)}(u)} + \Cov\sbra{ X^{(T)}(u), \sum_{i=N_R(uT)}^{N_R(vT)} W_i }
\label{}\\
=&\Var\sbra{X^{(T)}(u)}=\Var\sbra{\sum_{i=0}^{N_R(uT)}W_i},
\end{align}
which we can compute using the law of total variance:
\begin{align}
\Var\sbra{\frac{1}{\sqrt{T}}\sum_{i=0}^{N_R(uT)}W_i}=&\frac{1}{T}\bE\sbra{N_R(uT)}\Var\sbra{W}+\frac{1}{T}\Var\sbra{N_R(uT)}\rbra{\bE\sbra{W}}^2
\label{}\\
\end{align}
From Blackwell's renewal theorem (c.f. \cite[Sec. V, Theorem 4.4]{Asmussen2003}) we know that \linebreak $\frac{1}{T}\bE\sbra{N_R(uT)}\tend{}{T\rightarrow\infty}mu$. Thus,
\begin{align}
\Cov\sbra{X^{(T)}(u), X^{(T)}(v)}\tend{}{T\rightarrow\infty}u\rbra{\frac{m\alpha}{(1-\alpha)^3}+\frac{m^3\Var\sbra{\tau}}{(1-\alpha)^2}}=\sigma^2u.
\end{align}
}
Previously, we determined the convergence of the finite dimensional distributions and found their covariance structure. We can now present the proof of the central limit theorem.

\Proof[Proof of Theorem \ref{TheoremCLT}]{ In Lemma \ref{RHPLCLTSumprocess} we determined that the finite dimensional distributions are asymptotically normal. All there is left is to prove tightness. For this, let $\epsilon>0$. We want to show that
\begin{align}
\lim_{\delta\rightarrow 0}\limsup_{T\rightarrow\infty}\bP\rbra{\sup_{\absol{v-u}<\delta}\absol{X^{(T)}(v)-X^{(T)}(u)}>\epsilon}=0.
\end{align}
Without loss of generality, assume that $0<u<v<1$, then
\begin{align}
X^{(T)}(v)-X^{(T)}(u)=\frac{\sum_{i=N_R(uT)}^{N_R(vT)}W_i-\bE\sbra{\sum_{i=N_R(uT)}^{N_R(vT)}W_i}}{\sqrt{T}}.
\end{align}
Let us denote $\Delta N_R(T,u,v):=N_R(vT)-N_R(uT)$ to the number of clusters in $(uT,vT]$. Then, we compute the following variance
\begin{align}
\Var\sbra{\sum_{i=N_R(uT)}^{N_R(vT)}W_i}=\bE\sbra{\Delta N_R(T,u,v)}\frac{\alpha}{(1-\alpha)^3}+\Var\sbra{\Delta N_R(T,u,v)}\frac{1}{(1-\alpha)^2}.
\end{align}
Now, let $\rbra{\cF_t}_{t\ge 0}$ be the natural filtration of the RHP. Then $X^{(T)}(v)$ is adapted to the filtration  $\rbra{\cF_{Tv}}_{v\in[0,1]}$. Let $0<u<v<1$ and using the independence of the clusters that arrived after $uT$ from the clusters that arrived before $uT$, notice that
\begin{align}
\bE\sbra{X^{(T)}(v)-X^{(T)}(u)\mid \cF_{Tu}}=&\bE\sbra{\frac{\sum_{i=N_R(uT)}^{N_R(vT)}W_i-\bE\sbra{\sum_{i=N_R(uT)}^{N_R(vT)}W_i}}{\sqrt{T}}\Bigg| \cF_{Tu}}
\label{}\\
=&\frac{\bE\sbra{\sum_{i=N_R(uT)}^{N_R(vT)}W_i}-\bE\sbra{\sum_{i=N_R(uT)}^{N_R(vT)}W_i}}{\sqrt{T}}=0.
\end{align}
In other words,
\begin{align}
\bE\sbra{X^{(T)}(v)\mid \cF_{uT}}=X^{(T)}(u),
\end{align}
so $X^{(T)}(v)$ is a martingale adapted to $\rbra{\cF_{vT}}_{v\in[0,1]}$.
From Doob's martingale inequality, we have
\begin{align}
\bP\rbra{\supre{}{\absol{v-u}<\delta}\absol{X^{(T)}(v)-X^{(T)}(u)}>\epsilon}\le&\frac{\bE\sbra{\absol{X^{(T)}(v)-X^{(T)}(u)}^2}}{\epsilon^2}
\label{}\\
=&\frac{\Var\sbra{\sum_{i=N_R(uT)}^{N_R(vT)}W_i}}{T\epsilon^2}
\label{}\\
=&\frac{\bE\sbra{\Delta N_R(T,u,v)}\frac{\alpha}{(1-\alpha)^3}+\Var\sbra{\Delta N_R(T,u,v)}\frac{1}{(1-\alpha)^2}}{T\epsilon^2},
\end{align}
and from the elementary renewal theorem and central limit theorem for the renewal process, we obtain
\begin{align}
\frac{\bE\sbra{\Delta N_R(T,u,v)}\frac{\alpha}{(1-\alpha)^3}+\Var\sbra{\Delta N_R(T,u,v)}\frac{1}{(1-\alpha)^2}}{T\epsilon^2}\tend{}{T\rightarrow\infty}\frac{\frac{m\alpha(v-u)}{(1-\alpha)^3}+\frac{m^3(v-u)\Var\sbra{\tau}}{(1-\alpha)^2}}{\epsilon^2}.
\end{align}
This means that
\begin{align}
\limsup_{T\rightarrow\infty}\le \bP\rbra{\supre{}{\absol{v-u}<\delta}\absol{X^{(T)}(v)-X^{(T)}(u)}>\epsilon}\le\frac{\frac{m\alpha(v-u)}{(1-\alpha)^3}+\frac{m^3(v-u)\Var\sbra{\tau}}{(1-\alpha)^2}}{\epsilon^2},
\end{align}
and since the right hand side goes to zero as the difference $v-u$ becomes smaller, this entails that
\begin{align}
\lim_{\delta\rightarrow 0}\limsup_{T\rightarrow\infty} \bP\rbra{\supre{}{\absol{v-u}<\delta}\absol{X^{(T)}(v)-X^{(T)}(u)}>\epsilon}=0,
\end{align}
which proves tightness of the finite dimensional distributions, and therefore, the convergence 
\begin{align}
\rbra{\frac{1}{\sqrt{T}}\rbra{N(Tv)-\bE\sbra{N(Tv)}}}_{v\in[0,1]}\tend{d}{T\rightarrow\infty}\rbra{\sigma W(v)}_{v\in[0,1]},
\end{align}
in the Skorokhod topology.
}

Finally, we prove the asymptotic normality stated in Corollary \ref{CorollaryAsymptoticNormality}.
\Proof[Proof of Corollary \ref{CorollaryAsymptoticNormality}]{
Taking Lemma \ref{LemmaLLNSpeedConvergence} with $p=\frac{1}{2}$ yields
\begin{align}
\supre{}{v\in[0,1]}\absol{\frac{1}{\sqrt{T}}\bE\sbra{N(Tv)}-v\frac{m}{1-\alpha}\sqrt{T}}\tend{}{T\rightarrow\infty}0.
\end{align}
Furthermore, from Theorem \ref{TheoremCLT} we know that
\begin{align}
\rbra{\frac{1}{\sqrt{T}}N(Tv)-\frac{1}{\sqrt{T}}\bE\sbra{N(Tv)}}_{v\in[0,1]}\tend{d}{T\rightarrow\infty}\rbra{\sigma W(v)}_{v\in[0,1]},
\end{align}
from which the result follows.
}
Let us now use these results to obtain some insights about the behavior of the central limit theorem for the RHP in contrast to the classical Hawkes process. First of all, we notice that the variance consists of two parts
\begin{align}
\sigma^2=\frac{m\alpha}{(1-\alpha)^3}+\frac{m^3\Var\sbra{\tau}}{(1-\alpha)^2},
\end{align}
which have a simple interpretation. The first term corresponds to the variability of the points within a cluster, while the second term corresponds to the variability due to the arrivals of immigrants. We can see that if the renewal process is chosen with exponentially distributed inter arrivals, then $\Var\sbra{\tau}=m^{-2}$, which reduces to
\begin{align}
\sigma^2=m\rbra{\frac{\alpha+1-\alpha}{(1-\alpha)^3}}=\frac{m}{(1-\alpha)^3},
\end{align}
which corresponds to the variance of the classical Hawkes process. While, if we set $\alpha=0$, then
\begin{align}
\sigma^2=m^3\Var\sbra{\tau},
\end{align}
which is the variance for a renewal process in the central limit theorem.
We can see the effects on the variance of the RHP with a particular example. Let us consider a function $h$ that satisfies assumptions \textbf{(A0)} and \textbf{(A1)}. Let us consider an RHP with Weibull immigration. The interarrival distribution has the density function
\begin{align}
f(x)=bkx^{k-1}e^{-bx^k},\quad x\ge0,
\end{align}
where $b$ is called the scale parameter and $k$ the shape parameter. In this case, the hazard function for the renewals $\tau$ takes a simple polynomial form:
\begin{align}
\mu(x)=bkx^{k-1},\quad x\ge0.
\end{align}
The $n$-th raw moments of $\tau$ are given as
\begin{align}
\bE\sbra{\tau^n}=b^{-\frac{n}{k}}\Gamma\rbra{1+\frac{n}{k}},
\end{align}
and the variance is given as
\begin{align}
\Var\sbra{\tau}=b^{-\frac{2}{k}}\sbra{\Gamma\rbra{1+\frac{2}{k}}-\rbra{\Gamma\rbra{1+\frac{1}{k}}}^2}.
\end{align}
This variance becomes smaller as the shape parameter increases as seen below.
\begin{figure}[H]
    \centering
    \includegraphics[width=0.8\textwidth]{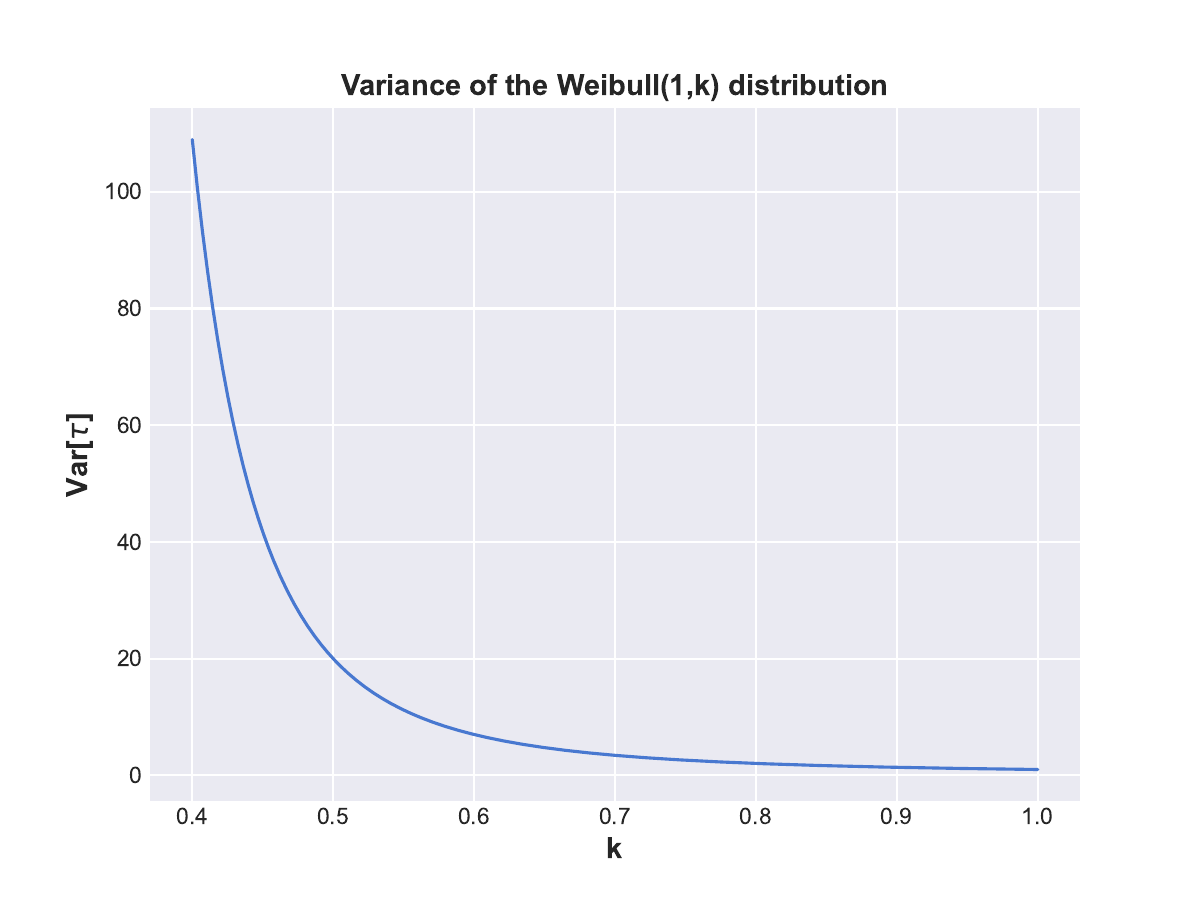}
    \caption{Variance of a Weibull random variable as a function of the shape parameter $k$ for fixed scale parameter $1$.}
    \label{fig:function_plot}
\end{figure}
Now, let us consider a family of RHP $\cbra{N_k}_{k\ge0}$ with Weibull immigration, but for each value of the shape parameter $k$, we set
\begin{align}
b(k)=\Gamma\rbra{1+\frac{1}{k}}^k,
\end{align}
so that in all cases $\bE\sbra{\tau}=1$. Then all these RHP have the same mean, but different values of $\sigma_R^2$. Notice that the case $k=1$ corresponds to the classical Hawkes process and it indeed recovers the variance for that case. As we increase the shape parameter of the Weibull distribution, the occurrence of clusters becomes more regular, so the variability due to the immigration decreases, while it increases as we make $k$ smaller. 

We can simulate renewal Hawkes processes using the R package `RHawkes' developed by Chen--Stindl based on the paper \cite{ChenDirect}. We include the case for $\alpha=0.5$ and Weibull(3,2) interarrivals.

\begin{figure}[H]
    \centering
    \includegraphics[width=0.8\textwidth]{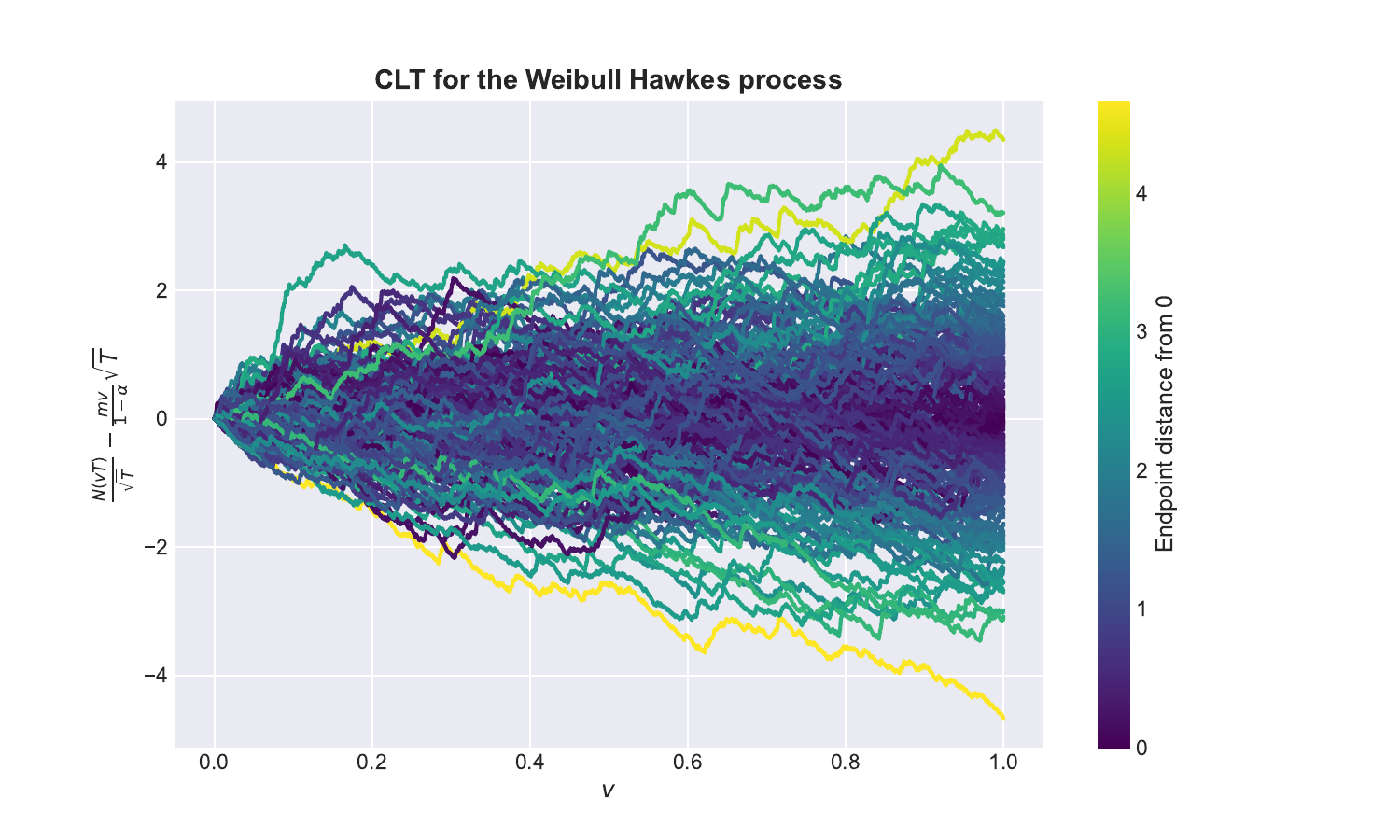}
    \caption{Limiting Brownian motions for an RHP with Weibull(3,2) inter arrivals.}
    \label{fig:CLT_plot}
\end{figure}
To find the variances, we can compute the sample variance at each time step and do a linear regression forced to pass through the origin. The slope will give us an approximation of $\sigma^2$.
\begin{figure}[H]
    \centering
    \includegraphics[width=0.8\textwidth]{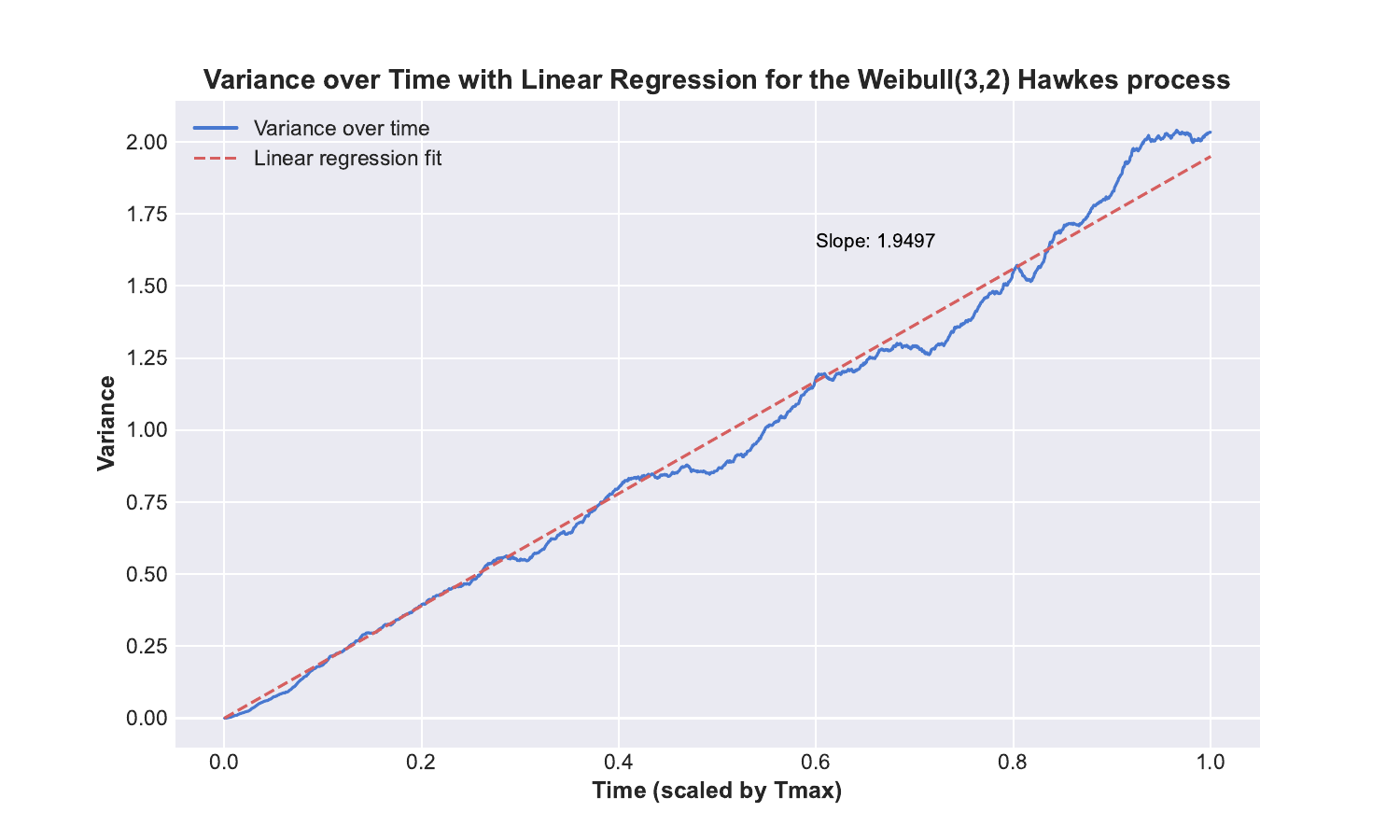}
    \caption{Numerical estimation of the variance for the Weibull(3,2) RHP.}
    \label{fig:Var_plot}
\end{figure}
By using the value of $\sigma^2$ obtained in Theorem \ref{TheoremCLT}, we have:
\begin{align}
\frac{0.5}{(1-0.5)^3}\times\frac{1}{3\Gamma(3/2)}+\frac{1}{(1-\alpha)^2}9\sbra{\Gamma\rbra{2}-\rbra{\Gamma\rbra{3/2}^2}}\approx 1.916.
\end{align}

\section{Acknowledgements}
I want to thank my PhD supervisor Dr. Kouji Yano for his guidance and comments on all the iterations of the manuscript for this paper. I want to thank as well professor Atsushi Takeuchi for his feedback on my work. I also extend my thanks to Dr. David Croydon, whose insightful comments helped me improve my results. I want to thanks as well the referees of this paper whose comments made it possible to improve the quality of the paper.

\bibliographystyle{plainreversed}

\end{document}